\documentclass[11pt]{article}
\usepackage{mathrsfs}
\usepackage{amssymb}
\usepackage[usenames,dvipsnames]{xcolor}
\usepackage{amsthm}
\usepackage{amsmath}
\usepackage{amsfonts}

\usepackage{graphicx}

 \numberwithin{equation}{section}

\def\bd{\mathbf{d}}

\newtheorem{theorem}{Theorem}[section]
\newtheorem{lemma}[theorem]{Lemma}
\newtheorem{corollary}[theorem]{Corollary}
\newtheorem{proposition}[theorem]{Proposition}

\theoremstyle{definition}
\newtheorem{definition}[theorem]{Definition}

\theoremstyle{remark}
\newtheorem{remark}[theorem]{Remark}

\begin{document}
\title{Space-time decay estimates for the incompressible viscous resistive Hall-MHD equations}
\author{Shangkun WENG\footnote{PDE and Functional Analysis Research Center (PARC), Seoul National University, 151-742, Seoul, Republic of Korea. {\it Email: skwengmath@gmail.com.}}}
\date{Seoul National University}
\maketitle

\def\be{\begin{eqnarray}}
\def\ee{\end{eqnarray}}
\def\p{\partial}
\def\no{\nonumber}
\def\e{\epsilon}
\def\de{\delta}
\def\De{\Delta}
\def\om{\omega}
\def\Om{\Omega}
\def\f{\frac}
\def\th{\theta}
\def\la{\lambda}
\def\b{\bigg}
\def\al{\alpha}
\def\La{\Lambda}
\def\ga{\gamma}
\def\Ga{\Gamma}
\def\ti{\tilde}
\def\Th{\Theta}
\def\si{\sigma}
\def\Si{\Sigma}
\def\bt{\begin{theorem}}
\def\et{\end{theorem}}
\def\bc{\begin{corollary}}
\def\ec{\end{corollary}}
\def\bl{\begin{lemma}}
\def\el{\end{lemma}}
\def\bp{\begin{proposition}}
\def\ep{\end{proposition}}
\def\br{\begin{remark}}
\def\er{\end{remark}}
\def\bd{\begin{definition}}
\def\ed{\end{definition}}
\def\bpf{\begin{proof}}
\def\epf{\end{proof}}
\begin{abstract}
  In this paper, we address the space-time decay properties for strong solutions to the incompressible viscous resistive Hall-MHD equations. We obtained the same space-time decay rates as those of the heat equation. Based on the temporal decay results in \cite{cs}, we find that one can obtain weighted estimates of the magnetic field $B$ by direct weighted energy estimate, and then by regarding the magnetic convection term as a forcing term in the velocity equations, we can obtain the weighted estimates for the vorticity, which yields the corresponding estimates for the velocity field. The higher order derivative estimates will be obtained by using a parabolic interpolation inequality proved in \cite{k01}. It should be emphasized that the the magnetic field has stronger decay properties than the velocity field in the sense that there is no restriction on the exponent of the weight. The same arguments also yield the sharp space-time decay rates for strong solutions to the usual MHD equations.
\end{abstract}

\begin{center}
\begin{minipage}{5.5in}

\

\end{minipage}
\end{center}

\section{Introduction and Main results}

In this paper we address the space time decay properties for strong solutions to the incompressible viscous resistive Hall-Magnetohydrodynamic equations.  The incompressible viscous resistive Hall-MHD equations take the following form:
\be\label{hmhd}\begin{cases}
\p_t u+u\cdot\nabla u+\nabla \pi=B\cdot\nabla B +\Delta u,\\
\p_t B-\nabla\times (u\times B)+\nabla\times ((\nabla\times B)\times B)=\Delta B,\\
\text{div} u=\text{div} B=0,
\end{cases}
\ee
where $u(x,t)=(u_1(x,t), u_2(x,t), u_3(x,t))$ and $B(x,t)=(b_1(x,t), b_2(x,t),b_3(x,t))$, $(x,t)\in \Bbb R^3\times[0, \infty)$, are the fluid velocity and magnetic field, $\pi=p+\f{1}{2} |B|^2$, where $p$ is the pressure. We will consider the Cauchy problem for (\ref{hmhd}), so we prescribe the initial data
\be\no
u(x,0)=u_0(x),\quad B(x,0)=B_0(x).
\ee
The initial data $u_0$ and $B_0$  satisfy the divergence free condition,
\be\no
\text{div}\, u_0(x)=\text{div} \,B_0(x)=0.
\ee
Recently, there are many researches on the Hall-MHD equations, concerning global weak solutions \cite{adfl,cdl}, local and global (small) strong solutions \cite{cdl,cwu,cl13,fhn}, and singularity formation in Hall-MHD \cite{cweng}. The application of Hall-MHD equations is mainly from the understanding of magnetic reconnection phenomena \cite{hg,lighthill}, where the topology structure of the magnetic field changes dramatically and the Hall effect must be included to get a correct description of this physical process. The Hall-MHD equations are also derived from a two-fluids Euler-Maxwell system for electrons and ions, through a set of scaling limits, see \cite{adfl}. They also provided a kinetic formulation for the Hall-MHD.

We will address the space-time decay properties of the strong solutions to (\ref{hmhd}). Let us briefly review the history of the study of the space-time decay properties of the solutions to the incompressible Navier-Stokes equations. Leray \cite{leray} proposed the problem whether the weak solution to the Navier-Stokes decay to zero in the $L^2$ norm. This was first solved positively by Kato \cite{kato} in the 2-D case and Scheonbek \cite{s1} in 3-D case. See also \cite{b04a,bv07a,gw02a,gw02b,km,ms,sw,wiegner} and the reference therein for more details. The Fourier splitting method introduced by Scheonbek was able to give the explicit decay rate, i.e. $\|u(\cdot,t)\|_{L^2}= O(t^{-n/2(1/p-1/2)})$ if $u_0\in L^p(\mathbb{R}^n)\cap L^2(\mathbb{R}^n) (1\leq p<2)$. The temporal decay rates for higher order norms of solutions to the incompressible Navier-Stokes was also investigated in \cite{sw}. Takahashi\cite{takahashi} first studied the space-time decay of solutions to the Navier-Stokes with nonzero forcing and zero initial data. Then the nonzero initial data case was addressed by \cite{ss00,agss}, however, their decay results are different from those of the heat equation. Kukavica \cite{k01} first obtained the sharp decay rate for the solution to the Navier-Stokes equation based on a parabolic interpolation inequality. This was further improved in Kukavica and Torres \cite{kt06}. More precisely, they established the sharp rates of decay for any weighted norm of higher order, i.e.
\be\no
\||x|^{r} D_x^b u(\cdot,t)\|_{L^2} =O(t^{-\gamma_0+\f r2-\f b2}),
\ee
for all $0\leq r\leq a$,where $D_x^b$ denote all the derivatives of order $b$, under the assumptions
\be\no
&&\|u(\cdot,t)\|_{L^2}= O(t^{-\gamma_0})\\\no
&&\||x|^a u(\cdot,t)\|_{L^2} = O(t^{-\gamma_0+\f a 2})\quad \text{for some $a>0$.}
\ee
The last assumption was also verified in their subsequent papers \cite{kt07,k09}.

Chae and Schonbek \cite{cs} investigated the temporal decay estimates for weak solutions to Hall-MHD system. They also obtained algebraic decay rates for higher order Sobolev norms of strong solutions to (\ref{hmhd}) with small initial data. It turned out that the Hall term does not affect the time asymptotic behavior, and the time decay rates behaved like those of the corresponding heat equation. Here we record their main results in the following.
\bt\label{cs1}({\bf Theorem 1.1 and 1.2 in \cite{cs}})
{\it
Let $(u_0, B_0)\in (L^1(\mathbb{R}^3)\cap L^2(\mathbb{R}^3))$ with $\text{div u}_0=\text{div B}_0=0$. Then there exists a weak solution $(u,B)$ to (\ref{hmhd}), which satisfies
\be\label{cs101}
\|u(t)\|_{L^2}+ \|B(t)\|_{L^2} \leq C_0 (t+1)^{-\f{3}{4}}.
\ee
If in addition, $(u_0,B_0)\in H^m(\mathbb{R}^3)$ for $m\in \mathbb{N}$ and $m\geq 3$ and $\|u_0\|_{H^m}+ \|B_0\|_{H^m}\leq K_1$ for some small constant $K_1$, then the solution $(u,B)$ will become strong and belong to $L^{\infty}(\mathbb{R}_+; H^m(\mathbb{R}^3))$ and also satisfy
\be\label{cs102}
\|D^m u(t)\|_{L^2}+ \|D^m B(t)\|_{L^2}\leq C_m (t+1)^{-\f{m}{2}-\f{3}{4}}
\ee
for all $t\geq T_*$. Here $C_m$ depends on $m$ and $C_0$.
}
\et

In this paper, we will address the spatial decay properties of the above strong solution $(u,B)$ in Theorem \ref{cs1}. For simplicity, we assume that the initial data $(u_0,B_0)$ belong to the Schwartz class $\mathcal{S}$, so that for any $a\geq 0$ and $b\in \mathbb{N}_0=\{0,1,2,\cdots\}$,
\be\no
\||x|^a D^b u_0\|_{L^2}<\infty,\quad \||x|^a D^b B_0\|_{L^2}<\infty.
\ee
Hence the solution $(u,B)$ will satisfy
\be\label{main00}
\|u(\cdot,t)\|_{L^2}+ \|B(\cdot,t)\|_{L^2} =O(t^{-\gamma_0})
\ee
and
\be\label{main01}
\|D^b u(\cdot,t)\|_{L^p}+ \|D^b B(\cdot,t)\|_{L^p} =O(t^{-\gamma_0-\f b2-\f 34(1-\f 2p)})
\ee
for any $2\leq p\leq \infty$ and $b\in \mathbb{N}_0$, where $\gamma_0=\f 34$. Note that (\ref{main01}) can be easily obtained from (\ref{cs102}) by interpolation.

Our main results is stated as follows.
\bt\label{main}
{\it
Let $(u,B)$ be the strong solution to (\ref{hmhd}) in Theorem \ref{cs1} with the initial data $(u_0,B_0)$ belong to the Schwartz class $\mathcal{S}$. Then we have the following weighted estimates for $u$ and $B$:
\be\label{main1}
\||x|^a D^b u(\cdot,t)\|_{L^p} =O(t^{-\gamma_0+\f{a}{2}-\f{b}{2}-\f{3}{4}(1-\f{2}{p})})
\ee
for any $b\in \mathbb{N}_0$ and $0\leq a <b+\f{5}{2}$ and $2\leq p\leq \infty$;
\be\label{main2}
\||x|^a D^b B(\cdot,t)\|_{L^p} =O(t^{-\gamma_0+\f a2-\f b2-\f 34 (1-\f 2p)})
\ee
for all $b\in \mathbb{N}_0$ and $a\geq 0$ and $2\leq p\leq \infty$. Furthermore, for the vorticity $\omega(t,x)=\text{curl u}(t,x)$, we have
\be\label{main3}
\||x|^a D^b \omega(\cdot,t)\|_{L^p} =O(t^{-\gamma_0+\f a2-\f b2-\f 12-\f 34 (1-\f 2p)})
\ee
for all $b\in \mathbb{N}_0$ and $a\geq 0$ and $2\leq p\leq \infty$.
}
\et

\br\label{main4}
{\it
We find that the spatial decay property of the magnetic field is stronger than that of the velocity field in the sense that there is no restriction on the exponent of the weight. This is basically due to the pressure term in the velocity equations. Note that the spatial decay of the voricity field is also much stronger than the velocity field.
}
\er

\br\label{main5}
{\it
One can relax the conditions on the initial data, i.e. there are constants $r>0$ and $k\in \mathbb{N}_0$, such that for all $0\leq a\leq r$ and $0\leq b\leq k$
\be\no
\||x|^a D^b u_0(\cdot)\|_{L^2}<\infty, \quad \||x|^a D^b B_0(\cdot)\|_{L^2} <\infty.
\ee
Then the conclusions in Theorem \ref{main} also hold with some obvious modification.
}
\er

Our proof basically follows the ideas introduced by Kukavica and Torres in a series of paper \cite{k01,kt06,kt07,k09}. However, we have some interesting new observations. First, we observe that one can directly obtain the weighted estimate of $B$ under the assumptions (\ref{main00}) and (\ref{main01}). Then we estimate the weighted norm of the vorticity $\omega=\text{curl} u$ by regarding $\text{curl $(B\cdot\nabla B)$}$ as a forcing term. It turns out that we can get the decay rate for the vorticity without any restriction on the exponent of the weight. By using the relation between $u$ and $\omega$, we can estimate the weighted norm of $u$ as in \cite{kt07,k09}. With the help of these weighted norm estimates of $u$ and $B$, we can use a parabolic type interpolation inequality proved in \cite{k01} to get the decay rates for the weighted norms of higher order derivatives. As above, we first estimate the magnetic field $B$ and then the velocity field regarding the magnetic convection term $B\cdot\nabla B$ as a forcing term. Here the existence of the Hall term requires a separate treatment of $\||x|^{a}\nabla B\|_{L^2}$, which will be used in the induction for higher order derivatives. Although the Hall term contains the second order derivative, it is a quadratic term, one can put the weight on another $B$, this is the reason why the Hall term does not affect the decay rates.

Our arguments certainly work in the usual MHD case, yielding the sharp space-time decay rates for strong solutions to the incompressible MHD. There are many previous studies on the time asymptotic behaviors for the solutions to the usual MHD equations \cite{as,hh,hx01,hx05b,sss,zhou}. In \cite{bv07b}, the authors studied the local in time persistence of space decay rates for the incompressible MHD, showing that if the initial magnetic field decays sufficiently fast, then the space decay rates of MHD solutions behave as that of Navier-Stokes solutions. On the other hand, if the initial magnetic field is poorly localized, then the magnetic field will govern the decay. Here our estimates also cover the space-time decay of higher order derivatives and the long time behavior, this result seems to be new, so we also include this as a theorem here.

\bt\label{main4}
{\it
Let $(u,B)$ is a strong solution to the incompressible viscous resistive MHD equations with initial data $(u_0,B_0)\in \mathcal{S}$ with $\text{div} u_0=\text{div} B_0=0$. Assume that (\ref{main00}) and (\ref{main01}) holds. Then we have the following weighted estimates for $u$ and $B$:
\be\label{main5}
\||x|^a D^b u(\cdot,t)\|_{L^p} =O(t^{-\gamma_0+\f{a}{2}-\f{b}{2}-\f{3}{4}(1-\f{2}{p})})
\ee
for any $b\in \mathbb{N}_0$ and $0\leq a <b+\f{5}{2}$ and $2\leq p\leq \infty$;
\be\label{main6}
\||x|^a D^b B(\cdot,t)\|_{L^p} =O(t^{-\gamma_0+\f a2-\f b2-\f 34 (1-\f 2p)})
\ee
for all $b\in \mathbb{N}_0$ and $a\geq 0$ and $2\leq p\leq \infty$. Furthermore, for the vorticity $\omega(t,x)=\text{curl} u(t,x)$, we have
\be\label{main7}
\||x|^a D^b \omega(\cdot,t)\|_{L^p} =O(t^{-\gamma_0+\f a2-\f b2-\f 12-\f 34 (1-\f 2p)})
\ee
for all $b\in \mathbb{N}_0$ and $a\geq 0$ and $2\leq p\leq \infty$.
}
\et

\br\label{main8}
{\it
It is well-known that the assumptions (\ref{main00}) and (\ref{main01}) will be satisfied either for global strong solutions to MHD with small data or global weak solutions to MHD equations after a large finite time. See \cite{as,cs,sss} for more details.
}
\er

The paper will be organized as follows. In section \ref{preliminary}, we present some lemmas which are needed in the weighted norm estimates. The weighted norm estimates for the solutions and higher order norms will be treated in section \ref{zero} and \ref{high} separately.

\section{Preliminary}\label{preliminary}

The following lemma is needed in the weighted estimates.
\bl\label{hall4}
{\it
Let $\alpha_0>1,\alpha_1<1, \alpha_2<1$ and $\beta_1, \beta_2<1$. Assume that a continuously differentiable function $F:[1,\infty)\rightarrow [0,\infty)$ satisfies
\be\no
F'(t) &\leq& C_0 t^{-\alpha_0} F(t) + C_1 t^{-\alpha_1} F(t)^{\beta_1}+ C_2 t^{-\alpha_2} F(t)^{\beta_2}+C_3 t^{\gamma_2-1},\quad t\geq 1\\\no
F(1)  &\leq& K_0
\ee
where $C_0,C_1,C_2,C_3, K_0\geq 0$ and $\gamma_i=\f{1-\alpha_i}{1-\beta_i}>0$ for $i=1,2$. Assume that $\gamma_1\geq \gamma_2$, then there exists a constant $C^*$ depending on $\alpha_0,\alpha_1,\beta_1,\alpha_2, \beta_2,K_0,C_i, i=1,\cdots,4$, such that $F(t)\leq C^* t^{\gamma_1}$ for $t\geq 1$.
}
\el
\bpf
This lemma is a simple variant of Lemma 2.2 in \cite{kt06}. For the reader's convenience, we present the proof following the idea in \cite{kt06}.
Let $t_0\geq 1$ be $C_0 t_0^{-(\alpha_0-1)}=\f{\gamma_1}{2}$. By Young's inequality, we have
\be\no
&&F'(t)\leq 2C_0 F(t)+ C_1 t^{-\f{\alpha_1}{1-\beta_1}}+ C_2 t^{-\f{\alpha_2}{1-\beta_2}}+ C_3 t^{\gamma_1-1},\quad t\geq 1,\\\no
&&F(1)\leq K_0.
\ee
Hence by the standard Gronwall's inequality, we know $F(t_0)\leq K_1$, where $K_1$ depending on $\alpha_0,\alpha_1,\beta_1,\alpha_2, \beta_2,K_0,C_i, i=1,\cdots,4$. Let $K>0$ be such that
\be\no
K\geq \max\b\{(C_1 2^{3+\beta_1}\gamma_1^{-1})^{\f{1}{1-\beta_1}},(C_2 2^{3+\beta_2}\gamma_1^{-1})^{\f{1}{1-\beta_2}}, K_1,\f{8C_0}{\gamma_1}\b\}.
\ee
Denote $R=\{t\geq t_0: F(t)\leq 2 K t^{\gamma_1}\}$. Since $F(t_0)\leq K_1\leq K$, $t_0\in R$, and by continuity, there exists a maximal interval $[t_0,b)\subset R$. We show that $b=\infty$.

Suppose $t_0<b<\infty$. Then $F(b)=2K b^{\gamma_1}$ and $F'(b)\geq G'(b)$, where $G(t)= 2K t^{\gamma_1}$. Note that
\begin{equation}\no
\begin{aligned}
G'(b) &\leq F'(b)\leq C_0 b^{-\alpha_0} 2K b^{\gamma_1} + C_1 b^{-\alpha_1} (2Kb^{\gamma_1})^{\beta_1}+ C_2 b^{-\alpha_2}(2K b^{\gamma_1})^{\beta_2} + C_3 b^{\gamma_2-1}\\\no
&\leq K\gamma_1 b^{\gamma_1-1}\b(2 C_0 \gamma_1^{-1} b^{1-\alpha_0}+ C_1\gamma_1^{-1} 2^{\beta_1} K^{\beta_1-1} b^{\gamma_1 \beta_1-\alpha_1-\gamma_1+1}\\\no
&\quad+C_1\gamma_1^{-1} 2^{\beta_2} K^{\beta_2-1} b^{\gamma_1 \beta_2-\alpha_2-\gamma_1+1} +\f{C_3}{K\gamma_1}\b)\\\no
&\leq K\gamma b^{\gamma_1-1}\b(1+C_1\gamma_1^{-1} 2^{\beta_1} K^{\beta_1-1}+C_2\gamma_1^{-1} 2^{\beta_2} K^{\beta_2-1}+ \f{C_3}{K\gamma_1}\b)\\\no
&\leq \f{3}{2} K\gamma_1 b^{\gamma_1-1}.
\end{aligned}
\end{equation}
However, $G'(b)=2 K\gamma_1 b^{\gamma_1-1}$, which is a contradiction. Hence we finish the proof.
\epf

The following two lemmas are needed in the weighted estimates for higher order norms. Both of them have been proved in \cite{k01} and \cite{kt06}. Here we omit the proof.

\bl\label{hall10}
{\it
Let $p\in [1,\infty]$ and $T>0$. Assume that $u\in L^{\infty}((0,T); L^p(\mathbb{R}^n))$ and $t(u_t-\Delta u)\in L^{\infty}((0,T); L^p(\mathbb{R}^n))$. Then $t^{1/2} \nabla u\in L^{\infty}((0,T); L^p(\mathbb{R}^n))$ and the inequality
\be\no
\sup_{t/2\leq\tau\leq t}\|\nabla u(\cdot,t)\|_{L^p}^2 &\leq& C(\sup_{t/4\leq \tau\leq t}\|u(\cdot,\tau)\|_{L^p})(\sup_{t/4\leq\tau\leq t}\|(u_t-\Delta u)(\cdot,\tau)\|_{L^p})\\\no
&\quad&+ \f{C}{t} \sup_{t/4\leq\tau\leq t} \|u(\cdot,\tau)\|_{L^p}^2
\ee
holds for every $t\in (0,T)$.
}
\el

\bl\label{hall11}
{\it
Let $\tau_0>0$ and assume that $F: [\tau_0,\infty)\rightarrow [0,\infty)$ satisfies $\sup_{\tau_0\leq \tau\leq A} F(\tau)<\infty$ for all $A>\tau_0$. If there exist $C_0>0$ and $\gamma\in \mathbb{R}$ such that
\be\label{hall111}
\sup_{t/2\leq \tau\leq t} F(\tau)^2 \leq C_0 t^{-2\gamma} + C_0 t^{-\gamma} \sup_{t/4\leq \tau\leq t} F(\tau),\quad t\geq 4\tau_0
\ee
then $F(t)=O(t^{-\gamma})$ as $t\rightarrow \infty$.
}
\el

\section{The weighted estimates for $u$ and $B$} \label{zero}

By Theorem \ref{cs1}, we may assume there exists a constant $\gamma_0=\f{3}{4}$, such that
\be\label{hall1}
\|u(\cdot,t)\|_{L^2}+ \|B(\cdot,t)\|_{L^2}= O(t^{-\gamma_0})\quad \text{as $t\rightarrow \infty$}.
\ee
Then by Gagliardo-Nirenberg inequality and (\ref{cs102}), we have for $2\leq p\leq \infty$
\be\label{hall2}
\|\p_{\alpha} u(\cdot,t)\|_{L^p}+ \|\p_{\alpha}B(\cdot,t)\|_{L^p}= O(t^{-\gamma_0-\f{|\alpha|}{2}-\f{3}{4}(1-\f{2}{p})}),\quad \alpha\in \mathbb{N}_0^3.
\ee

First we observe that the weighted estimate for the magnetic field can be obtained directly under the assumptions (\ref{hall1}) and (\ref{hall2}).
\bt\label{hall3}
{\it
Under the assumption (\ref{hall1}) and (\ref{hall2}), we have
\be\label{hall301}
\||x|^a B(\cdot,t)\|_{L^2}= O(t^{-\gamma_0+a/2})\quad \text{as $t\rightarrow \infty$}
\ee
for all $a\geq 0$.
}
\et

\bpf
\be\no
B_t+ u\cdot\nabla B+ \nabla\times ((\nabla\times B)\times B)= B\cdot\nabla u+ \Delta B.
\ee
Multiplying the above equations by $2|x|^{2a} B$, and setting $G(t)= \int_{\mathbb{R}^3} |x|^{2a} |B(x,t)|^2 dx$, then we get
\be\no
&\quad&\f{d}{dt} G(t)+ 2\int_{\mathbb{R}^3} |x|^{2a} |\nabla B|^2dx \\\no
&=&- \int_{\mathbb{R}^3} 2 |x|^{2a} B \cdot(u\cdot\nabla B) dx + \int_{\mathbb{R}^3} 2 |x|^{2a} B\cdot (B\cdot\nabla u) dx\\\no
&\quad&-\int_{\mathbb{R}^3}2 |x|^{2a} B\cdot [\nabla\times((\nabla \times B)\times B)] dx- \int_{\mathbb{R}^3} 4a |x|^{2a-2} \sum_{i,j=1}^n B_i x_j \p_j B_i dx\\\no
&:=& I+ II+ III+ IV.
\ee
Then we estimate these fourth terms as follows.
\be\no
|I|&\leq& 2\int_{\mathbb{R}^3} |x|^{2a} |\nabla B| |B| |u| dx \leq \f{1}{3} \int_{\mathbb{R}^3} |x|^{2a} |\nabla B|^2 dx+C \int_{\mathbb{R}^3} |x|^{2a} |B|^2 |u|^2 dx \\\no
&\leq& \f{1}{3} \int_{\mathbb{R}^3} |x|^{2a} |\nabla B|^2 dx+ C\|u\|_{L^{\infty}}^2 G(t),
\ee
\be\no
|II|&\leq& 2\int_{\mathbb{R}^3} |x|^{2a} |\nabla u| |B|^2 dx\leq 2\|\nabla u\|_{L^{\infty}} G(t),
\ee
\be\no
|III|&\leq& C\int_{\mathbb{R}^3} |x|^{2a} |B| (|\nabla^2 B| |B|+ |\nabla B|^2) dx \\\no
&\leq& C\|\nabla^2 B\|_{L^{\infty}} G(t)+ C\|\nabla B\|_{L^{\infty}} \int_{\mathbb{R}^3} |x|^{2a} |B||\nabla B| dx \\\no
&\leq& \f{1}{3}\int_{\mathbb{R}^3} |x|^{2a} |\nabla B|^2 dx+C (\|\nabla^2 B\|_{L^{\infty}}+\|\nabla B\|_{L^{\infty}}^2) G(t).
\ee
\be\no
|IV|&\leq& C\int_{\mathbb{R}^3} |x|^{2a-1} |B| |\nabla B| dx \leq \f{1}{3} \int_{\mathbb{R}^3} |x|^{2a} |\nabla B|^2 dx + C \int_{\mathbb{R}^3}|x|^{2a-2}|B|^2 dx\\\no
&\leq& \f{1}{3} \int_{\mathbb{R}^3} |x|^{2a}|\nabla B|^2 dx+ C G(t)^{\f{a-1}{a}} \|B\|_{L^2}^{\f{2}{a}}.
\ee
Combining all these estimates, we get
\be\no
G'(t)&\leq& C(\|u\|_{L^{\infty}}^2+ \|\nabla u\|_{L^{\infty}}+ \|\nabla^2 B\|_{L^{\infty}}+ \|\nabla B\|_{L^{\infty}}^2)G(t)+C G(t)^{\f{a-1}{a}} \|B\|_{L^2}^{\f{2}{a}}\\\no
&\leq& C t^{-\gamma_0-\f{1}{2}-\f 34} G(t)+ C t^{-\f{2\gamma_0}{a}} G(t)^{\f{a-1}{a}}.
\ee
If $a> 2\gamma_0$, then we may apply Lemma \ref{hall4} with $\alpha_0=\gamma_0+\f{1}{2}+\f 34>1, \alpha_1=\f{2\gamma_0}{a}<1, \beta_1=\f{a-1}{a}<1, C_2=C_3=0$ to get $G(t)\leq C t^{\gamma_1}$ with
\be\no
\gamma_1= \f{1-\alpha_1}{1-\beta_1}== a-2\gamma_0
\ee
and the theorem is proved for all $a>2\gamma_0$. The conclusion for $a\in (0,2\gamma_0]$ follows by interpolation.
\epf

Now we turn to the velocity field. Let the vorticity $\omega(t,x)=\text{curl} u(t,x)$, then
\be\no
\p_t \omega + u\cdot \nabla \omega- \omega\cdot\nabla u-\Delta \omega= \text{curl}(B\cdot \nabla B).
\ee
With the weighted estimates (\ref{hall301}) of $B$ at hand, we regard $\text{curl}(B\cdot\nabla B)$ as a forcing term, and estimate the weighted norm of the vorticity as above.

\bt\label{hall15}
{\it
Under the assumptions (\ref{hall1}) and (\ref{hall2}), we have the following estimate for all $a\geq 0$
\be\label{hall151}
\||x|^a\omega(\cdot,t)\|_{L^2}= O(t^{-\gamma_0-\f 12+ \f a2}).
\ee
}
\et

\bpf
Multiplying the vorticity equation by $2 |x|^{2a} \omega$ and setting $F(t)= \int_{\mathbb{R}^3} |x|^{2a} |\omega(x,t)|^2 dx$, then we get
\be\no
&\quad&\f{d}{dt} F(t) + 2\int_{\mathbb{R}^3} |x|^{2a} |\nabla \omega(x,t)|^2 dx \\\no
&=&- \int_{\mathbb{R}^3} 2|x|^{2a} \omega\cdot(u\cdot\nabla\omega) dx+ \int_{\mathbb{R}^3} 2|x|^{2a}\omega\cdot (\omega\cdot\nabla u) dx\\\no
&\quad&-4a \int_{\mathbb{R}^3} |x|^{2a-2} \sum_{i,j=1}^n x_j \omega_i \p_j\omega_i dx+ \int_{\mathbb{R}^3} 2|x|^{2a} \omega \cdot \text{curl}(B\cdot\nabla B) dx\\\no
&:=& I+ II + III+ IV.
\ee
These four terms will be estimated as follows.
\be\no
|I|&\leq& \f{1}{3}\int_{\mathbb{R}^3} |x|^{2a} |\nabla \omega|^2 dx+ C\|u\|_{L^{\infty}}^2 F(t),
\ee
\be\no
|II|&\leq& 2\|\nabla u\|_{L^{\infty}} F(t),
\ee
\be\no
|III|&\leq& C \int_{\mathbb{R}^3} |x|^{2a-1} |\omega| |\nabla \omega| dx \\\no
&\leq& \f{1}{3}\int_{\mathbb{R}^3} |x|^{2a} |\nabla \omega|^2 dx+ C F(t)^{\f{a-1}{a}} \|\omega\|_{L^2}^{\f{2}{a}},
\ee
\be\no
|IV|&=&2\b|\int_{\mathbb{R}^3} \text{curl}(|x|^{2a}\omega)\cdot (B\cdot\nabla B) dx\b|\\\no
&=& 2\b|\int_{\mathbb{R}^3} |x|^{2a} \text{curl} \omega \cdot (B\cdot\nabla B) dx+ \int_{\mathbb{R}^3}2a|x|^{2a-2}(x\times \omega)\cdot (B\cdot\nabla B) dx\b|\\\no
&\leq& C\int_{\mathbb{R}^3} |x|^{2a} |\nabla \omega||B||\nabla B| dx + C \int_{\mathbb{R}^3} |x|^{2a-1} |\omega| |B||\nabla B| dx\\\no
&\leq& \f{1}{3} \int_{\mathbb{R}^3} |x|^{2a} |\nabla\omega|^2 dx+ C\|\nabla B\|_{L^{\infty}}^2 \int_{\mathbb{R}^3} |x|^{2a} |B|^2 dx+  F(t)^{\f{1}{2}}\|\nabla B\|_{L^{\infty}} \||x|^{a-1}|B|\|_{L^2}.
\ee
Combining all these estimates together, we obtain
\be\no
F'(t)&\leq& C_0(\|u\|_{L^{\infty}}^2+\|\nabla u\|_{L^{\infty}})F(t)+ C_0 F(t)^{\f{a-1}{a}} \|\omega\|_{L^2}^{\f{2}{a}} \\\no
&\quad&\quad\quad + C_0 F(t)^{\f{1}{2}} \|\nabla B\|_{L^{\infty}} \||x|^{a-1} B\|_{L^2}+ C_0\|\nabla B\|_{L^{\infty}}^2 \||x|^a B\|_{L^2}^2 \\\no
&\leq& C_0 t^{-2\gamma_0-5/4}F(t)+ C_0 F(t)^{\f{a-1}{a}} t^{-\f{2}{a}(\gamma_0+\f{1}{2})}+C_0 F(t)^{\f{1}{2}} t^{-2\gamma_0-\f{7}{4}+\f{a}{2}}+ C_0 t^{-4\gamma_0-\f{5}{2}+a}\\\no
&=& C_0t^{-2} F(t)+ C_0 F(t)^{\f{a-1}{a}} t^{-\f{5}{2a}}+ C_0 F(t)^{\f{1}{2}} t^{-\f{13}{4}+\f{a}{2}} + C_0 t^{a-\f{11}{2}}.
\ee
Now we can apply Lemma \ref{hall4}. Here $\beta_1=\f{a-1}{a}, \alpha_1=\f{5}{2a}, \beta_2=\f{1}{2}, \alpha_2=-\f{a}{2}+\f{13}{4}$. To assure that $\alpha_1<1,\alpha_2<1$, we require $a>\f{9}{2}$. Hence $\gamma_1=\f{1-\alpha_1}{1-\beta_1}=a-\f{5}{2}>\gamma_2=\f{1-\alpha_2}{1-\beta_2}=a-\f{9}{2}$. By Lemma \ref{hall4}, we obtain
\be\no
F(t)\leq C t^{-\f 52+a}=C t^{-2\gamma_0-1+a}.
\ee
\epf

By the relation $-\Delta u=\text{curl } \omega$ and the Caffarelli-Kohn-Nirenberg inequality \cite{ckn84}
\be\no
\||x|^a u\|_{L^p} \leq \||x|^{1+a} \nabla u\|_{L^p},
\ee
one can argue as in \cite{kt07} and \cite{k09} to obtain the weighted estimates for the velocity field $u$ as stated in the following theorem. Since the proof are almost the same, here we omit the details.
\bt\label{hall16}
{\it
Under the assumptions (\ref{hall1})-(\ref{hall2}), we have the following weighted estimates
\be\label{hall161}
\||x|^a u(\cdot,t)\|_{L^2}= O(t^{-\gamma_0+\f{a}{2}})
\ee
for all $a\in [0, \f{5}{2})$.
}
\et

\section{The Weighted estimates for higher order derivatives}\label{high}

Based on the estimates (\ref{hall301}) and (\ref{hall161}), we can apply Lemma \ref{hall10} to get the weighted estimates for higher order derivatives of $u$ and $B$. First we estimate $B$, here the existence of the hall term requires a separate treatment of $\||x|^{a}\nabla B\|_{L^2}$, which will be used for the induction of higher order derivatives. Although the Hall term contains the second order derivative, it is a quadratic term, one can put the weight on another $B$, this is the reason why the Hall term does not affect the decay rates.

\bt\label{hall21}
{\it
Under the assumptions  (\ref{hall1}) and (\ref{hall2}), then the following estimates hold for all $a\geq 0, b\in \mathbb{N}_0$ and $2\leq p\leq\infty$
\be\label{hall210}
\||x|^a D^b B(\cdot,t)\|_{L^p}= O(t^{-\gamma_0-\f{b}{2}+\f{a}{2}-\f{3}{4}(1-\f{2}{p})}).
\ee
}
\et

\bpf
We only need to prove the case $a>2$, since we already know (\ref{hall210}) holds for $a=0$, the case $0<a\leq 2$ can be obtained by interpolation. By the Gagliardo-Nirenberg inequality, the case $p>2$ follows from $p=2$. Indeed, for any $f\in L^2(\mathbb{R}^3)\cap \dot{H}^2(\mathbb{R}^3)$, one has $\|f\|_{L^{\infty}(\mathbb{R}^3)}\leq \|f\|_{L^2(\mathbb{R}^3)}^{\f 14}\|f\|_{\dot{H}^2(\mathbb{R}^3)}^{\f 34}$. Hence one can derive the estimate of $\||x|^a D^b B(\cdot,t)\|_{L^{\infty}(\mathbb{R}^3)}$ from those of $\||x|^{a} D^b B(\cdot,t)\|_{L^2(\mathbb{R}^3)}$. The case $p\in (2,\infty)$ just follows from interpolation. Therefore, we assume $p=2$.
For $a>2$, we choose the weight $\phi$:
\be\no
\phi(x,t)= (|x|^2+t)^{\f{a}{2}},\quad t\geq 1,
\ee
then by simple calculations, we get
\be\no
|\nabla\phi(x,t)|\leq (|x|^2+t)^{\f{a-1}{2}},\quad |(\p_t-\Delta)\phi(x,t)|\leq C(|x|^2+t)^{\f{a}{2}-1}.
\ee
The case $b=0$ has been proved in Theorem \ref{hall3}. The following proof will be separated into two steps. The first step addresses the case $b=1$, which is needed in the weighted estimate for higher order derivatives in the second step.

{\it Step 1.} The case $b=1$. To use Lemma \ref{hall10}, we need to derive the equation for $\phi B$
\be\no
(\p_t-\Delta)(\phi B)= (\p_t-\Delta)\phi B-2\nabla \phi\cdot \nabla B-\phi(u\cdot\nabla B-B\cdot\nabla u+ \nabla \times((\nabla\times B)\times B)).
\ee
Hence we obtain
\begin{equation}\no
\begin{aligned}
\sup_{t/2\leq \tau\leq t}\|\nabla(\phi B)\|_{L^2}^2 &\leq C\sup_{t/4\leq \tau\leq t}\|\phi B\|_{L^2} \sup_{t/4\leq \tau\leq t}\b(\|(\p_t\phi-\Delta\phi) B\|_{L^2}+ \|\nabla \phi\cdot \nabla B\|_{L^2}\\\no
&\quad+ \|\phi u\cdot\nabla B\|_{L^2}+ \|\phi \nabla\times((\nabla\times B)\times B)\|_{L^2}+ \|\phi B\cdot\nabla u\|_{L^2}\b)\\\no
&\quad+\f{C}{t}\sup_{t/4\leq \tau\leq t} \|\phi B(\tau)\|_{L^2}^2
\end{aligned}
\end{equation}
Note that
\be\no
\sup_{t/2\leq \tau\leq t}\|\nabla(\phi B)\|_{L^2}^2&\geq& \f{1}{2}\sup_{t/2\leq \tau\leq t}\|\phi\nabla B\|_{L^2}^2-C\sup_{t/2\leq \tau\leq t}\|\nabla \phi B\|_{L^2}^2\\\no
&\geq& \f{1}{2}\sup_{t/2\leq \tau\leq t}\|\phi\nabla B\|_{L^2}^2-O(t^{-2\gamma_0+a-1}),\\\no
\|(\p_t-\Delta)\phi B\|_{L^2}&\leq& \|(|x|^2+t)^{\f{a}{2}-1}B\|_{L^2}\leq\||x|^{a-2}B\|_{L^2}+t^{\f{a}{2}-1}\|B\|_{L^2}\\\no
&\leq& O(t^{-\gamma_0+\f{a}{2}-1}),
\ee
\be\no
\|\nabla \phi\cdot\nabla B\|_{L^2}&\leq&O(t^{-\f{1}{2}})\|\phi \nabla B\|_{L^2},\\\no
\|\phi u\cdot\nabla B\|_{L^2}&\leq& \|u\|_{L^{\infty}}\|\phi\nabla B\|_{L^2}\leq O(t^{-\gamma_0-\f{3}{4}})\|\phi \nabla B\|_{L^2},\\\no
\|\phi B\cdot\nabla u\|_{L^2}&\leq& \|\phi B\|_{L^2} \|\nabla u\|_{L^{\infty}}\leq O(t^{-2\gamma_0+\f{a}{2}-\f{3}{4}-\f{1}{2}}).
\ee
\be\no
\|\phi \nabla\times((\nabla\times B)\times B)\|_{L^2}&\leq& \|\phi |\nabla^2 B||B|\|_{L^2}+\|\phi |\nabla B|^2\|_{L^2}\\\no
&\leq& \|\nabla^2B\|_{L^{\infty}}\|\phi B\|_{L^2}+ \|\nabla B\|_{L^{\infty}}\|\phi \nabla B\|_{L^2}\\\no
&\leq& O(t^{-2\gamma_0+\f{a}{2}-\f{3}{4}-1})+ O(t^{-\gamma_0-\f{1}{2}-\f{3}{4}})\|\phi\nabla B\|_{L^2}.
\ee
hence we obtain
\begin{equation}
\begin{aligned}
\sup_{t/2\leq\tau\leq t}\|\phi\nabla B\|_{L^2}^2 &\leq O(t^{-2\gamma_0+a-1})+O(t^{-\gamma_0+\f{a}{2}})\b(O(t^{-\gamma_0+\f{a}{2}-1})+O(t^{-2\gamma_0+\f{a}{2}-\f{3}{4}-1})\\\no
&+(O(t^{-\f12})+O(t^{-\gamma_0-\f34}))\sup_{t/4\leq\tau\leq t}\|\phi \nabla B\|_{L^2}\b)+\f{C}{t}O(t^{-2\gamma_0+a})\\\no
&\leq O(t^{-2\gamma_0+a-1})+ O(t^{-\gamma_0+\f{a}{2}-\f{1}{2}})\sup_{t/4\leq\tau\leq t}\|\phi \nabla B\|_{L^2}.
\end{aligned}
\end{equation}
By applying Lemma \ref{hall11}, we get
\be\no
\||x|^a \nabla B(\cdot,t)\|_{L^2}= O(t^{-\gamma_0-\f{1}{2}+\f{a}{2}}).
\ee
{\it Step 2.} The case $b\geq 2$. Assume that the conclusion holds for all the derivatives up to order $b\geq 1$, we want to show that it also holds for $b+1$. Take any $\alpha\in \mathbb{N}_0^3$ with $|\alpha|=b$, then
\be\no
(\p_t-\Delta)(\phi \p_{\alpha} B)&=& (\p_t\phi-\Delta \phi)\p_{\alpha}B- 2\p_j\phi \p_j \p_{\alpha} B\\\no
&\quad&-\sum_{0\leq \beta\leq \alpha} C_{\alpha,\beta} \phi \b(\p_{\beta} u_j\p_j \p_{\alpha-\beta} B-\p_{\beta} B_j\p_j\p_{\alpha-\beta} u\b)\\\no
&\quad&-\sum_{0\leq \beta\leq \alpha}C_{\alpha,\beta} \nabla\times \b(\p_{\beta}(\nabla\times B)\times (\p_{\alpha-\beta} B)\b).
\ee
Hence by Lemma \ref{hall10}, we obtain
\begin{equation}\no
\begin{aligned}
\sup_{t/2\leq \tau\leq t}\|\nabla(\phi \p_{\alpha} B)\|_{L^2}^2 &\leq C\sup_{t/4\leq \tau\leq t}\|\phi \p_{\alpha} B\|_{L^2} \sup_{t/4\leq \tau\leq t}\b[\|(\p_t\phi-\Delta\phi) \p_{\alpha}B\|_{L^2}+ \|\nabla \phi\cdot \nabla \p_{\alpha}B\|_{L^2}\\\no
&\quad+ \sum_{0\leq \beta\leq \alpha}C_{\alpha,\beta}\b(\|\phi \p_{\beta}u_j\p_j\p_{\alpha-\beta}B\|_{L^2}+ \|\phi \p_{\beta}B_j\p_j\p_{\alpha-\beta}u\|_{L^2}\\\no
&\quad+\|\phi \nabla\times((\nabla\times \p_{\beta}B)\times \p_{\alpha-\beta}B)\|_{L^2}\b)\b]+\f{C}{t}\sup_{t/4\leq \tau\leq t} \|\phi \p_{\alpha}B(\tau)\|_{L^2}^2
\end{aligned}
\end{equation}
By induction assumptions, we have
\be\no
\sup_{t/2\leq \tau\leq t}\|\nabla(\phi \p_{\alpha} B)(\tau)\|_{L^2}^2&\geq&\f{1}{2}\sup_{t/2\leq \tau\leq t}\|\phi\nabla\p_{\alpha} B\|_{L^2}^2-O(t^{-2\gamma_0-b+a-1}),\\\no
\sup_{t/4\leq \tau\leq t}\|\phi \p_{\alpha} B\|_{L^2}&\leq& O(t^{-\gamma_0-\f{b}{2}+\f{a}{2}}),\quad\sup_{t/4\leq \tau\leq t}\|(\p_t-\Delta)\phi \p_{\alpha} B\|_{L^2}\leq O(t^{-\gamma_0-\f b2+\f a2-1}),\\\no
\sup_{t/4\leq \tau\leq t}\|\nabla \phi\nabla\p_{\alpha} B\|_{L^2}&\leq& O(t^{-\f 12}) \sup_{t/4\leq \tau\leq t}\|\phi\nabla \p_{\alpha} B\|_{L^2}.
\ee
For the other three terms, we estimate as follows
\begin{equation}\no
\begin{aligned}
\|\phi\p_{\beta}u_j\p_j \p_{\alpha-\beta} B\|_{L^2} &\leq \begin{cases}
\|\p_{\beta} u_j\|_{L^{\infty}} \|\phi \p_j\p_{\alpha-\beta} B\|_{L^2},\quad \text{if $|\beta|>0$}\\
\|u\|_{L^{\infty}} \|\phi \nabla \p_{\alpha} B\|_{L^2},\quad \text{if $\beta=0$}
\end{cases}\\\no
&\leq\begin{cases}
O(t^{-2\gamma_0-\f b2+\f a2-\f 12-\f 34}),\quad \text{if $|\beta|>0$}\\
O(t^{-\gamma_0-\f 34})\|\phi \nabla \p_{\alpha} B\|_{L^2} ,\quad \text{if $\beta=0$}
\end{cases}\\\no
\|\phi\p_{\beta} B_j \p_j\p_{\alpha-\beta}u\|_{L^2}&\leq \|\phi\p_{\beta} B_j\|_{L^2}\|\p_j\p_{\alpha-\beta}u\|_{L^{\infty}} \\\no
&\leq O(t^{-2\gamma_0-\f b2+\f a2-\f 34-\f 12}).
\end{aligned}
\end{equation}
\be\no
\|\phi \nabla\times((\nabla\times \p_{\beta}B)\times \p_{\alpha-\beta}B)\|_{L^2}&\leq&\|\phi |\p_{\beta}\nabla^2 B||\p_{\alpha-\beta}B|\|_{L^2}+ \|\phi|\p_{\beta}\nabla B||\p_{\alpha-\beta} \nabla B|\|_{L^2}\\\no
&:=& J_1+J_2.
\ee
\begin{equation}\no
\begin{aligned}
J_1&\leq \begin{cases}
\|\phi\p_{\beta}\nabla^2 B\|_{L^{2}} \|\p_{\alpha-\beta} B\|_{L^{\infty}},\quad \text{if $|\beta|\leq b-2$}\\
\|\p_{\beta} \nabla^2 B\|_{L^{\infty}} \|\phi \p_{\alpha-\beta} B\|_{L^2},\quad \text{if $|\beta|=b-1$ or $b$}
\end{cases}\\\no
&\leq O(t^{-2\gamma_0-\f b2+\f a2-\f 34-1}),\\\no
J_2&\leq \begin{cases}
\|\phi\p_{\beta}\nabla B\|_{L^{2}} \|\p_{\alpha-\beta}\nabla B\|_{L^{\infty}},\quad \text{if $|\beta|\leq b-1$}\\
\|\p_{\alpha} \nabla B\|_{L^{\infty}} \|\phi \nabla B\|_{L^2},\quad \text{if $\beta=\alpha$}
\end{cases}\\\no
&\leq O(t^{-2\gamma_0-\f b2+\f a2-\f 34-1}).
\end{aligned}
\end{equation}
Note that in the estimate of $J_2$, we have used the result from {\it Step 1}.

Combining all these estimates together, we get
\be\no
\sup_{t/2\leq\tau\leq t}\|\phi \nabla \p_{\alpha}B\|_{L^2}^2 &\leq& O(t^{-2\gamma_0-b+a-1})+ O(t^{-\gamma_0-\f b2+\f a2-\f 12})\sup_{t/4\leq\tau\leq t}\|\phi \nabla \p_{\alpha}B\|_{L^2}^2.
\ee
This implies
\be\no
\|\phi \nabla \p_{\alpha} B\|_{L^2}= O(t^{-\gamma_0-\f b2+\f a2-\f 12}).
\ee
\epf

\bt\label{hall28}
{\it
Under the assumptions (\ref{hall1}) and (\ref{hall2}), we have
\be\no
\||x|^a D^b u(\cdot,t)\|_{L^p} =O(t^{-\gamma_0-\f{b}{2}-\f{3}{4}(1-\f{2}{p})+\f{a}{2}})
\ee
for any $a\in [0,\f{5}{2})$ and $b\in \mathbb{N}_0$ and $p\in [2,\infty]$.
}
\et

\bpf
As before, we only consider $a\geq 2$. By the Gagliardo-Nirenberg inequality, the case $p>2$ follows from $p=2$. Therefore, we assume $p=2$. Denote
\be\no
\phi(x,t)= (|x|^2+t)^{\f{a}{2}}.
\ee
The case $b=0$ has been proved in Lemma \ref{hall16}. Assume that the conclusion holds for $b\in \mathbb{N}_+$, and we shall establish it for $b+1$. Fix $\alpha\in \mathbb{N}_0^3$ such that $|\alpha|=b$. Then we have
\be\no
(\p_t-\Delta)(\phi \p_{\alpha} u)&=& (\p_t\phi-\Delta \phi)\p_{\alpha}u- 2\p_j\phi \p_j \p_{\alpha}u- \phi \nabla \p_{\alpha} \pi\\\no
&\quad&-\sum_{0\leq \beta\leq \alpha} C_{\alpha,\beta} \phi \b(\p_{\beta} u_j\p_j \p_{\alpha-\beta} u-\p_{\beta} B_j\p_j\p_{\alpha-\beta} B\b),
\ee
Then Lemma \ref{hall10} yields
\begin{equation}\no
\begin{aligned}
\sup_{t/2\leq \tau\leq t} \|\nabla (\phi \p_{\alpha}u)\|_{L^2}^2 &\leq \sup_{t/4\leq\tau\leq t}\|\phi \p_{\alpha}u\|_{L^2} \sup_{t/4\leq \tau\leq t}\b(\|(\p_t-\Delta)\phi \p_{\alpha}u\|_{L^2}+ \|\nabla \phi\cdot\nabla \p_{\alpha}u\|_{L^2}\\\no
&\quad+\|\phi\nabla \p_{\alpha}\pi\|_{L^2}+\sum_{0\leq \beta\leq \alpha}C(\alpha,\beta)(\|\phi\p_{\beta}u_j\p_j\p_{\alpha-\beta}u\|_{L^2}\\\no
&\quad+\|\phi\p_{\beta}B_j\p_j\p_{\alpha-\beta}B\|_{L^2})\b)+\f{C}{t} \sup_{t/4\leq\tau\leq t}\|\phi \p_{\alpha} u\|_{L^2}^2.
\end{aligned}
\end{equation}
As above, we have
\be\no
\sup_{t/2\leq \tau\leq t}\|\nabla(\phi \p_{\alpha} u)(\tau)\|_{L^2}^2&\geq&\f{1}{2}\sup_{t/2\leq \tau\leq t}\|\phi\nabla\p_{\alpha} u\|_{L^2}^2-O(t^{-2\gamma_0-b+a-1}),\\\no
\sup_{t/4\leq \tau\leq t}\|\phi \p_{\alpha} u\|_{L^2}&\leq& O(t^{-\gamma_0-\f{b}{2}+\f{a}{2}}),\quad\sup_{t/4\leq \tau\leq t}\|(\p_t-\Delta)\phi \p_{\alpha} u\|_{L^2}\leq O(t^{-\gamma_0-\f b2+\f a2-1}),\\\no
\sup_{t/4\leq \tau\leq t}\|\nabla \phi\nabla\p_{\alpha} u\|_{L^2}&\leq& O(t^{-\f 12}) \sup_{t/4\leq \tau\leq t}\|\phi\nabla \p_{\alpha} u\|_{L^2}.
\ee
And
\begin{equation}\no
\begin{aligned}
\|\phi\p_{\beta}u_j\p_j \p_{\alpha-\beta} u\|_{L^2} &\leq \begin{cases}
\|\p_{\beta} u_j\|_{L^{\infty}} \|\phi \p_j\p_{\alpha-\beta} u\|_{L^2},\quad \text{if $|\beta|>0$}\\
\|u\|_{L^{\infty}} \|\phi \nabla \p_{\alpha} u\|_{L^2},\quad \text{if $\beta=0$}
\end{cases}\\\no
&\leq\begin{cases}
O(t^{-2\gamma_0-\f b2+\f a2-\f 12-\f 34}),\quad \text{if $|\beta|>0$}\\
O(t^{-\gamma_0-\f 34})\|\phi \nabla \p_{\alpha} u\|_{L^2} ,\quad \text{if $\beta=0$}
\end{cases}\\\no
\|\phi\p_{\beta} B_j \p_j\p_{\alpha-\beta}B\|_{L^2}&\leq \|\phi\p_{\beta} B_j\|_{L^2}\|\p_j\p_{\alpha-\beta}B\|_{L^{\infty}} \\\no
&\leq O(t^{-2\gamma_0-\f b2+\f a2-\f 34-\f 12}).
\end{aligned}
\end{equation}
For $\|\phi \nabla\p_{\alpha}\pi\|_{L^2}$, note that
\be\no
\pi = R_i R_j(u_i u_j-B_i B_j),
\ee
we need to use the weighted estimates for Riesz operator in Lemma 4.2 in \cite{k01}. Hence
\be\no
\|\phi \nabla\p_{\alpha}\pi\|_{L^2} &\leq& \||x|^a \nabla \p_{\alpha}\pi\|_{L^2}+ t^{\f a2}\|\nabla \p_{\alpha}\pi\|_{L^2} \\\no
&\leq& \||x|^a\nabla \p_{\alpha}(u_i u_j-B_iB_j)\|_{L^2}+ \||x|^{a-1}\p_{\alpha}(u_i u_j-B_i B_j)\|_{L^2}+ t^{\f a2}\|\nabla \p_{\alpha}\pi\|_{L^2}\\\no
&:=& H_1+H_2+H_3,
\ee
\be\no
H_1&\leq& \sum_{0\leq \beta\leq \alpha}\b(\|\phi |\nabla \p_{\beta} u||\p_{\alpha-\beta}u|\|_{L^2}+ \|\phi |\nabla \p_{\beta} B||\p_{\alpha-\beta}B|\|_{L^2}\b)\\\no
&\leq& \|\phi u\|_{L^2} \|\nabla \p_{\alpha} u\|_{L^{\infty}}+\|\phi B\|_{L^2} \|\nabla \p_{\alpha} B\|_{L^{\infty}}\\\no
&\quad&+ \sum_{0\leq \beta<\alpha}\b(\|\phi\nabla \p_{\beta}u\|_{L^2}\|\p_{\alpha-\beta}u\|_{L^{\infty}}+\|\phi\nabla \p_{\beta}B\|_{L^2}\|\p_{\alpha-\beta}B\|_{L^{\infty}}\b)\\\no
&\leq& O(t^{-2\gamma_0-\f b2+\f a2-\f 34-\f 12}),\\\no
H_2&\leq& \sum_{0\leq \beta\leq \alpha} \b(\||x|^{a-1}\p_{\beta} u\|_{L^2}\|\p_{\alpha-\beta} u\|_{L^{\infty}}+\||x|^{a-1}\p_{\beta} B\|_{L^2}\|\p_{\alpha-\beta} B\|_{L^{\infty}}\b)\\\no
&\leq& O(t^{-2\gamma_0-\f b2+\f a2-\f 34-\f 12}),\\\no
H_3&\leq& t^{\f a2}\|\nabla\p_{\alpha}(u_i u_j-B_i B_j)\|_{L^2}\\\no
&\leq& t^{\f a2}\sum_{0\leq \beta\leq \alpha}\b(\|\nabla\p_{\beta} u\|_{L^2}\|\p_{\alpha-\beta} u\|_{L^{\infty}}+\|\nabla\p_{\beta} B\|_{L^2}\|\p_{\alpha-\beta} B\|_{L^{\infty}}\b)\\\no
&\leq& O(t^{-2\gamma_0-\f b2+\f a2-\f 34-\f 12}).
\ee
Hence we obtain
\be\no
\sup_{t/2\leq \tau\leq t} \|\phi \nabla \p_{\alpha} u\|_{L^2}^2 &\leq& O(t^{-2\gamma_0-b+a-1})+ O(t^{-\gamma_0-\f b2+\f a2-\f 12})\sup_{t/4\leq\tau\leq t}\|\phi \nabla \p_{\alpha}u\|_{L^2}^2.
\ee
Applying Lemma \ref{hall11}, we get
\be\no
\|\phi \nabla \p_{\alpha} u\|_{L^2}= O(t^{-\gamma_0-\f b2-\f 12+\f a2}).
\ee

\epf

Now we show that the vorticity field has much stronger decay properties than the velocity field in the sense that there is no restriction on the exponent of the weight. Indeed, we have
\bt\label{hall26}
{\it
Under the assumptions (\ref{hall1}) and (\ref{hall2}), the following estimates
\be\label{hall261}
\||x|^a D^b \omega(x,t)\|_{L^p}= O(t^{-\gamma_0-\f b2-\f 12+\f a2-\f{3}{4}(1-\f{2}{p})})
\ee
hold for any $a\geq 0$ and $b\in \mathbb{N}_0$ and $2\leq p\leq \infty$.
}
\et

\bpf
We choose same weight function as before. The conclusion is true for $b=0$ as showed in Theorem \ref{hall15}. We assume that the conclusion holds for any derivatives up to order $b$, we want to show that it also holds for $b+1$. Take any $\alpha\in \mathbb{N}_0^3$ with $|\alpha|=b$, then
\begin{equation}\no
\begin{aligned}
(\p_t-\Delta)(\phi \p_{\alpha}\omega)&=(\p_t-\Delta)\phi \p_{\alpha}\omega- 2(\nabla\phi\cdot\nabla)\p_{\alpha}\omega-\sum_{0\leq\beta\leq \alpha}C(\alpha,\beta)\b(\phi \p_{\beta} u_j\p_j\p_{\alpha-\beta}\omega\\\no
&\quad-\phi \p_{\beta}\omega_j\p_j\p_{\alpha-\beta}u-\phi \nabla\times(\p_{\beta}B_j \p_j\p_{\alpha-\beta} B)\b).
\end{aligned}
\end{equation}
Then applying Lemma \ref{hall10}, we obtain
\begin{equation}\no
\begin{aligned}
\sup_{t/2\leq \tau\leq t}\|\nabla(\phi \p_{\alpha} \omega)\|_{L^2}^2 &\leq C\sup_{t/4\leq \tau\leq t}\|\phi \p_{\alpha} \omega\|_{L^2} \sup_{t/4\leq \tau\leq t}\b[\|(\p_t\phi-\Delta\phi) \p_{\alpha}\omega\|_{L^2}+ \|\nabla \phi\cdot \nabla \p_{\alpha}\omega\|_{L^2}\\\no
&\quad+ \sum_{0\leq \beta\leq \alpha}C_{\alpha,\beta}\b(\|\phi \p_{\beta}u_j\p_j\p_{\alpha-\beta}\omega\|_{L^2}+ \|\phi \p_{\beta}\omega_j\p_j\p_{\alpha-\beta}u\|_{L^2}\\\no
&\quad+\|\phi \nabla\times(\p_{\beta}B_j\p_j\p_{\alpha-\beta}B)\|_{L^2}\b)\b]+\f{C}{t}\sup_{t/4\leq \tau\leq t} \|\phi \p_{\alpha}\omega\|_{L^2}^2
\end{aligned}
\end{equation}
As above, we have
\be\no
\sup_{t/2\leq \tau\leq t}\|\nabla(\phi \p_{\alpha} \omega)\|_{L^2}^2&\geq&\f{1}{2}\sup_{t/2\leq \tau\leq t}\|\phi\nabla\p_{\alpha} \omega\|_{L^2}^2-O(t^{-2\gamma_0-b-1+a-1}),\\\no
\sup_{t/4\leq \tau\leq t}\|\phi \p_{\alpha} \omega\|_{L^2}&\leq& O(t^{-\gamma_0-\f{b}{2}-\f 12+\f{a}{2}}),\quad\sup_{t/4\leq \tau\leq t}\|(\p_t-\Delta)\phi \p_{\alpha} \omega\|_{L^2}\leq O(t^{-\gamma_0-\f b2-\f 12+\f a2-1}),\\\no
\sup_{t/4\leq \tau\leq t}\|\nabla \phi\nabla\p_{\alpha} \omega\|_{L^2}&\leq& O(t^{-\f 12}) \sup_{t/4\leq \tau\leq t}\|\phi\nabla \p_{\alpha} \omega\|_{L^2}.
\ee
For the other three terms, we estimate as follows
\begin{equation}\no
\begin{aligned}
\|\phi\p_{\beta}u_j\p_j \p_{\alpha-\beta} \omega\|_{L^2} &\leq \begin{cases}
\|\p_{\beta} u_j\|_{L^{\infty}} \|\phi \p_j\p_{\alpha-\beta} \omega\|_{L^2},\quad \text{if $|\beta|>0$}\\
\|u\|_{L^{\infty}} \|\phi \nabla \p_{\alpha} \omega\|_{L^2},\quad \text{if $\beta=0$}
\end{cases}\\\no
&\leq\begin{cases}
O(t^{-2\gamma_0-\f b2+\f a2-1-\f 34}),\quad \text{if $|\beta|>0$}\\
O(t^{-\gamma_0-\f 34})\|\phi \nabla \p_{\alpha} \omega\|_{L^2} ,\quad \text{if $\beta=0$}
\end{cases}\\\no
\|\phi\p_{\beta} \omega_j \p_j\p_{\alpha-\beta}u\|_{L^2}&\leq \|\phi\p_{\beta} \omega_j\|_{L^2}\|\p_j\p_{\alpha-\beta}u\|_{L^{\infty}} \\\no
&\leq O(t^{-2\gamma_0-\f b2+\f a2-\f 34-1}).
\end{aligned}
\end{equation}
\be\no
\|\phi \nabla\times(\p_{\beta}B_j\p_j\p_{\alpha-\beta}B)\|_{L^2}&\leq&\|\phi |\p_{\beta}\nabla B||\p_j\p_{\alpha-\beta}B|\|_{L^2}+ \|\phi|\p_{\beta}B||\p_{\alpha-\beta} \nabla^2 B|\|_{L^2}\\\no
&\leq&\|\phi\nabla \p_{\beta}B_j\|_{L^2}\|\nabla\p_{\alpha-\beta}B\|_{L^{\infty}}+ \|\phi\p_{\beta}B\|_{L^2}\|\nabla^2\p_{\alpha-\beta}B\|_{L^{\infty}}\\\no
&\leq&O(t^{-2\gamma_0-\f b2+\f a2-\f 34-1}).
\ee
Combining all these estimates together, we get
\be\no
\sup_{t/2\leq\tau\leq t}\|\phi\nabla \p_{\alpha}\omega\|_{L^2}^2 &\leq& O(t^{-2\gamma_0-b+a-2})+ O(t^{-\gamma_0-\f b2+\f a2-1})\sup_{t/4\leq\tau\leq t}\|\phi\nabla \p_{\alpha}\omega\|_{L^2}.
\ee
Hence $\|\phi \p_{\alpha}\omega\|_{L^2}= O(t^{-\gamma_0-\f{b+1}{2}-\f 12+\f a2})$.
\epf

In particular, we have showed that $\||x|^a u(\cdot,t)\|_{L^2}=O(t^{-\gamma_0+\f a2})$ for all $a\in [0,\f 52)$. Now one can argue as in Theorem 3.2 of \cite{kt06} to get the following theorem.
\bt\label{hall35}
{\it
Assume that (\ref{hall1}) and (\ref{hall2}) hold, then
\be\label{hall351}
\||x|^a D_x^b u(\cdot,t)\|_{L^p} = O(t^{-\gamma_0-\f b2+\f a2-\f 34(1-\f 2p)})
\ee
for every $2\leq p\leq \infty$, $a \in[0, b+\f{5}{2})$ and $b\in \mathbb{N}_0$.
}
\et

Indeed, the proof of Theorem 3.2 in \cite{kt06} only produce the estimate for $p=2$. The other cases can be derived from $p=2$ by Gagliardo-Nirenberg inequality and interpolation theory.

{\it Proof of Theorem \ref{main} and \ref{main4}.} Theorem \ref{hall21}, \ref{hall26} and \ref{hall35} imply Theorem \ref{main}. From the proof of Theorem \ref{hall21}, \ref{hall26} and \ref{hall35}, one can see that same arguments certainly work for the usual incompressible viscous, resisitve MHD equations, hence Theorem \ref{main4} holds.

{\bf Acknowledgements.} The author would like to thank Prof. Dongho Chae for his interest and encouragement. Special thanks also go to Prof. Brandolese for his interest in this paper and kindly notifying the author of his important works on space-time decay for Navier-Stokes and MHD.

\bibliographystyle{plain}

\begin{thebibliography}{10}
\bibitem{adfl}
M. Acheritogaray, P. Degond, A. Frouvelle, J.-G. Liu. {\it  Kinetic formulation and global existence for the Hall-Magneto-hydrodynamics system,}
Kinet. Relat. Models 4 (2011) 901--918.

\bibitem{as}
R. Agapito, M. Schonbek. {\it Non-uniform decay of MHD equations with and without magnetic diffusion.} Comm. Partial Differential Equations 32 (2007), no. 10-12, 1791--1812.

\bibitem{agss}
C. Amrouche,V. Girault, M. Schonbek, T. Schonbek. {\it Pointwise decay of solutions and of higher derivatives to Navier-Stokes equations.} SIAM J. Math. Anal. 31 (2000), no. 4, 740-753.

\bibitem{b04a}
L. Brandolese. {\it Space-time decay of Navier-Stokes flows invariant under rotations.} Math. Ann. 329 (2004), no. 4, 685--706.

\bibitem{bv07a}
L. Brandolese, F. Vigneron. {\it New asymptotic profiles of nonstationary solutions of the Navier-Stokes system.} J. Math. Pures Appl. (9) 88 (2007), no. 1, 64--86.

\bibitem{bv07b}
L. Brandolese, F. Vigneron. {\it On the localization of the magnetic and the velocity fields in the equations of magnetohydrodynamics}. Proc. Roy. Soc. Edinburgh Sect. A 137 (2007), no. 3, 475--495.

\bibitem{ckn84}
L. Caffarelli, R, Kohn, L. Nirenberg. {\it First order interpolation inequalities with weights.} Compositio Math. 53 (1984), 259-275.

\bibitem{cl13}
D. Chae, J. Lee. {\it On the blow-up criterion and small data global existence for the Hall-magnetohydrodynamics.} Journal of Differential Equations, Volume 256, Issue 11, 2014, 3835--3858.

\bibitem{cs}
D. Chae, M. Schonbek. {\it On the temporal decay for the Hall-magnetohydrodynamic equations,} J. Differential Equations 255 (2013), no. 11, 3971--3982.

\bibitem{cdl}
D. Chae, P. Degond, J.-G. Liu. {\it  Well-posedness for Hall-magnetohydrodynamics,}  Ann. Inst. H. Poincar¨¦ Anal. Non Lin¨¦aire 31 (2014), no. 3, 555--565.

\bibitem{cweng}
D. Chae, S. Weng. {\it Singularity formation for the incompressible Hall-MHD equations without resistivity.} arXiv:1312.5519.

\bibitem{cwu}
D. Chae, J. Wu. {\it Local well-posedness for the Hall-MHD equations with fractional magnetic diffusion.}  arXiv:1404.0486.



\bibitem{fhn}
J. Fan, S. Huang, G. Nakamura. {\it  Well-posedness for the axisymmetric incompressible viscous Hall-magnetohydrodynamic equations,}  Appl. Math. Lett. 26 (2013), no. 9, 963--967.

\bibitem{gw02a}
T. Gallay, C. Wayne. {\it Invariant manifolds and the long-time asymptotics of the Navier-Stokes and vorticity equations on $\mathbb{R}^2$.} Arch. Ration. Mech. Anal. 163 (2002), no. 3, 209--258.

\bibitem{gw02b}
T. Gallay, C. Wayne. {\it Long-time asymptotics of the Navier-Stokes and vorticity equations on $\mathbb{R}^3$.} Recent developments in the mathematical theory of water waves (Oberwolfach, 2001). R. Soc. Lond. Philos. Trans. Ser. A Math. Phys. Eng. Sci. 360 (2002), no. 1799, 2155--2188.

\bibitem{hh}
P. Han, C. He. {\it Decay properties of solutions to the incompressible magnetohydrodynamics equations in a half space.} Math. Methods Appl. Sci. 35 (2012), no. 12, 1472--1488.

\bibitem{hx01}
C. He, Z. Xin. {\it On the decay properties of solutions to the non-stationary Navier-Stokes equations in $\mathbb{R}^3$.} Proc. Roy. Soc. Edinburgh Sect. A 131 (2001), no. 3, 597--619.
\bibitem{hx05b}

C. He, Z. Xin. {\it Partial regularity of suitable weak solutions to the incompressible magnetohydrodynamic equations.} J. Funct. Anal. 227 (2005), no. 1, 113--152.

\bibitem{hg}
H. Homann, R. Grauer. {\it Bifurcation analysis of magnetic reconnection in Hall-MHD
systems,} Physica D 208 (2005), pp. 59-72.

\bibitem{leray}
J. Leray. {\it Sur le mouvement d'un liquide visqueux emplissant l'espace.} (French) Acta Math. 63 (1934), no. 1, 193--248.

\bibitem{lighthill}
M. J. Lighthill, {\it Studies on magneto-hydrodynamic waves and other anisotropic wave motions,}  Philos. Trans. R. Soc. Lond. Ser. A 252 (1960)
397--430.

\bibitem{km}
R. Kajikiya, T. Miyakawa. {\it On $L^2$ decay of weak solutions of the Navier-Stokes equations in $\mathbb{R}^n$.} Math. Z. 192 (1986), no. 1, 135--148.

\bibitem{kato}
T. Kato. {\it Strong $L^p$-solutions of the Navier-Stokes equation in $\mathbb{R}^m$, with applications to weak solutions.} Math. Z. 187 (1984), 471-480.


\bibitem{k01}
I. Kukavica. {\it Space-time Decay for Solutions of the Navier-Stokes Equations.} Indiana University Mathematics Journal, vol. 50, no. 1(2001), 205-222.

\bibitem{kt06}
I. Kukavica; J. J. Torres. {\it Weighted bounds for the velocity and the vorticity for the Navier-Stokes equations.} Nonlinearity 19 (2006), no. 2, 293--303.

\bibitem{kt07}
I. Kukavica; J. J. Torres. {\it Weighted $L^p$ decay for solutions of the Navier-Stokes equations.} Comm. Partial Differential Equations 32 (2007), no. 4-6, 819--831.

\bibitem{k09}
I. Kukavica. {\it On the weighted decay for solutions of the Navier-Stokes system.} Nonlinear Anal. 70 (2009), no. 6, 2466--2470.

\bibitem{ms}
T. Miyakawa, M. Schonbek. {\it On optimal decay rates for weak solutions to the Navier-Stokes equations in $\mathbb{R}^n$.} Math. Bohem. 126 (2001), no. 2, 443--455.

\bibitem{s1}
M. Schonbek. {\it $L^2$ decay for weak solutions of the Navier-Stokes equations.} Arch. Rational Mech. Anal. 88 (1985), no. 3, 209-222.

\bibitem{sw}
M. Schonbek, M. Wiegner. {\it On the decay of higher-order norms of the solutions of Navier-Stokes equations.} Proc. Roy. Soc. Edinburgh Sect. A 126 (1996), no. 3, 677-685.

\bibitem{sss}
M. Schonbek, T. Schonbek, E. S$\ddot{u}$li. {\it Large-time behaviour of solutions to the magnetohydrodynamics equations.} Math. Ann. 304 (1996), no. 4, 717-756.

\bibitem{ss00}
M. Schonbek, T. Schonbek. {\it On the boundedness and decay of moments of solutions to the Navier-Stokes equations.} Adv. Differential Equations 5 (2000), no. 7-9, 861-898.


\bibitem{takahashi}
S. Takahashi. {\it A weighted equation approach to decay rate estimates for the Navier-Stokes equations.} Nonlinear Anal. 37 (1999), no. 6, Ser. A: Theory Methods, 751--789.

\bibitem{wiegner}
M. Wiegner. {\it Decay results for weak solutions of the Navier-Stokes equations on $\mathbb{R}^n$.} J. London Math. Soc. (2) 35 (1987), no. 2, 303--313.

\bibitem{zhou}
Y, Zhou. {\it Remarks on regularities for the 3D MHD equations.} Discrete Contin. Dyn. Syst. 12 (2005), no. 5, 881--886.
\end{thebibliography}

\end{document}